\documentclass[12pt]{article}
\usepackage[cp1251]{inputenc}
\usepackage[english]{babel}
\usepackage[T2B]{fontenc}

\usepackage[pdftex, colorlinks, urlcolor=blue, pdfborder={0 0 0 [2 0]}]{hyperref} % для arxiv.org

\usepackage{amsfonts,amsmath,amssymb}%,theorem
\usepackage{amsthm}
\usepackage{array}
\usepackage{esint}
\usepackage{graphicx}

\usepackage{misccorr} % в заголовках появляется точка, но при ссылке на них  её нет

\newtheorem{theorem}{\hspace{\parindent}\bf{Theorem}}%[section]

\newtheorem{m-lemma}{\hspace{\parindent}\bf{Микролемма}}%[section]

\newtheorem{dfn}{\hspace{\parindent}\sl{D\,e\,f\,i\,n\,i\,t\,i\,o\,n\,\,}}%[section]

\newenvironment{proofs}
{\vspace{1pt}\par{\sl%
P\,r\,o\,o\,f\.\,\ }}%
{\noindent\vspace{1pt}}

{\vspace{1pt}\par{\sl%
Д\,о\,к\,а\,з\,а\,т\,е\,л\,ь\,с\,т\,в\,о\ \ т\,е\,о\,р\,е\,м\,ы\, }}%
{\noindent\vspace{1pt}}

{\noindent\vspace{1pt}}

\newenvironment{proof-m}
{\vspace{1pt}\par{\sl%
Д\,о\,к\,а\,з\,а\,т\,е\,л\,ь\,с\,т\,в\,о\ \ м\,и\,к\,р\,о\,л\,е\,м\,м\,ы.\,\ }}%
{\noindent\vspace{1pt}}

{\vspace{1pt}\par{\sl%
Р\,е\,ш\,е\,н\,и\,е.\,}}%
{\noindent\vspace{1pt}}

\makeatletter
\renewcommand{\@biblabel}[1]{#1.} % Заменяем библиографию с квадратных скобок на точку:
\makeatother
\begin{document}

\title{ Stochastic First Integrals, Kernel Functions for Integral Invariants and the Kolmogorov equations}

\author{Valery Doobko,  Elena Karachanskaya  \\
Academy of Municipal Management, Kyiv\\
Pacific National University, Khabarobsk}

%\address{Department of foundational and Computer Sciences, Pacific National University, Russia, Khabarovsk}

\date{}

\maketitle
\begin{abstract}
In this article the authors present stochastic first integrals (SFI), the generalized It\^o-Wentzell formula and its application for obtaining  the equations for SFI, for  kernel functions for integral invariants and the Kolmogorov equations, described by the generalized It\^o equations.

\end{abstract}


\section*{Introduction}
The concept of а first integral for solution of the deterministic dynamical system  is the fundamental concept of an analytical mechanics. For stochastic differential equations (SDE) similar concepts exist as well. They are {\it a  first integral} for  It\^o's SDE (Doobko, 1978 \cite{D_78-e}); {\it a first   forward integral} and {\it a first   backward integral} for  It\^o's SDE (Krylov and Rozovsky, 1982 \cite{KR_82-e}); {\it a stochastic first   integral} for generalized It\^o's SDE (GSDE) \cite{D_02-e,11_KchIto-e}. However a direct transfer of the concept of а first integral from deterministic systems to stochastic systems is impossible.

The different conservation laws, such as  energy, weight, impulse, impulse moment, etc., are the basis of invariants and first integrals.  For example, if a collection of the enumerable number of the initial solutions for the same dynamical equation is connected to points which are similar to  particles, then the number of this points is a conservative value since the conditions of the existence and uniqueness of the solution are fulfilled.

The limiting state of this representation is a density of these points number and conserving of the integral for it into real space.  A function which has the same property is called a kernel function for an integral invariant.

By imposing certain restrictions a partial  differential equation for the kernel function can be obtained \cite{D_02-e,D_89-e,13_KchArxDir-e}. And we do consider a situation when properties of the  functionals which are conserved  on some boundary  space-time continuum  domain can be restrictions.
For such cases we refer to representation of local invariants  \cite{D_95-e}.  Evolving structures and  functionals connected with initial values domain, and are considered as {\it  dynamical invariants}. For example, an element of a phase-space volume and  hypersurfaces are dynamical invariants.

The theory of stochastic integral invariants (SIIs) is one of the  approaches for studying   stochastic dynamical systems described by the  It\^o equations \cite{D_84-e,D_89-e,D_98-e}. This approach is very effective for  defining the conserved functionals for  evolving  systems, such as: first integrals and stochastic first integrals for the  It\^o equations (SDE) \cite{D_78-e,D_02-Dan-e}, the length of the random chain  \cite{10_KchTSP}, a constant velocity of the random walking points   \cite{97_DubChalBroun-e}. It is applied for  constructing  the analytical solution of the Langevin type equation \cite{98_DubChalSolut-e}; for a determination of probability moments of the point  which is walking randomly  on a sphere  \cite{10_KchHf-e,12_KchSphera-e}; for  obtaining the It\^o-Wentzell formula for SDE  and the generalized  It\^o-Wentzell formula for generalized SDE  (GSDE ) \cite{D_89-e,D_02-e,11_KchOboz-e,12_KchDubPrep-e}; for forming and  solving the problem of the control program  with probability one (PCP1) for stochastic systems which are subjected to strong perturbations \cite{07_ChOboz-e,09_ChUprW-e,11_KchUpr-e}.

We will show  an application of the concept of the SFI  for forming  and proving the  theorems about forms of the equations for the SFI. Further, the existence and uniqueness theorems for solutions of GSDE for kernel functions of  integral invariants are defined. The important point of this research is the generalized It\^o-Wentzell formula. This formula is a differentiation law for compound random function which depends on solutions of  GSDE \cite{D_02-e,11_KchIto-e,12_KchDubPrep-e}.  The choice of such notation comes from the fact that without Poisson perturbations this formula rearranges to the well-known It\^o-Wentzell formula \cite{KR_82-e,D_89-e,Wentzel_65-e}.

The aim of this article is to demonstrate  a possibility of applying the generalized It\^o-Wentzell formula for obtaining equations for the SFI, for a complete proof of deriving the equation for stochastic kernel functions (SKFs) of the SIIs and the derivation of the Kolmogorov  equations \cite{13_DKchUch-e}.

The article includes three sections. In the first section the following results were obtained:  the concept of the SFI is define, the generalized It\^o-Wentzell formula is demonstrated  and equations for the SFI are given. In the second section  we consider a concept of a local stochastic density of a dynamical invariant which is connected with a solution of the GSDE, then we set an equation for it and we establish a link between  the local stochastic density and the concept of the SKF for a SII. In the third section we yield the Kolmogorov equation.

\section{Stochastic first integrals}

Let $(\Omega,\mathcal{F}, \mathbf{P})$ be a complete probability space, ${\rm {\mathcal F}}_{t}=\bigl\{{\mathcal F}_{t},  \ \ t\in [0,T]\bigr\}$,   ${\mathcal F}_{s}\subset {\mathcal F}_{t} $, $s<t$ be a non decreasing flow of the $\sigma-$algebras.

Let us consider the next random processes: ${\bf w}(t)$  is the $m$-dimensional Wiener process;
$\nu (\Delta t;\Delta \gamma )$ is the standard Poisson measure defined on the space $ [0;T]\times\mathbb{R}^{n'} $ and has properties:   ${\rm M}[\nu (\Delta t,\Delta \gamma )]=\Delta t\cdot \Pi (\Delta \gamma )$; $\int_{{\mathbb R} (\gamma )}  \Pi (d\gamma )<\infty $, $ {\mathbb R} (\gamma )=\mathbb{R}^{n'}$ (it means a finiteness of jumps intensity during a small time interval);
the one-dimensional processes $w_{k}(t)$ and the Poisson measure $\nu([0;T],\mathcal{A})$ are defined on the same space $(\Omega,\mathcal{F},\mathbf{P})$ and $\mathcal{F}_{t}-$measurable for every $t>0$ for each set $\mathcal{A}$ from $\sigma-$algebras of the Borel sets and are mutually independent.

Let a random process $x(t)$ defined on ${\mathbb R}^{n} $, be a solution of the SDE \cite{GS_68-e}:
\begin{equation} \label{Ayd01.1u}
\begin{array}{c}
\displaystyle{dx_{i} (t)=a_{i} (t )
dt+b_{i k} (t )
dw_{k} (t)+\int_{{\mathbb R} (\gamma )} g_{i} (t; \gamma ) \nu (dt;d\gamma ),} \\
\displaystyle {x(t)=x(t;x_{0})\bigr|_{t=0} =x_{0}, \ \ \ {\mbox {\rm для всех}} \ \ x_{0}\in {\rm {\mathbb R}}^{n} }, \ \ i=\overline{{1,n}};\ k=\overline{{1,m}},
\end{array}
\end{equation}
where
$i =\overline{{1,n}}$, $k=\overline{{1,m}} $, $ {\mathbb R} (\gamma )=\mathbb{R}^{n'}$ is a space of vectors  $\gamma$;
and summation is held on the indices double appeared.

The coefficients  $a(t;x)$, $b(t;x)$ and $g(t;x;\gamma )$ satisfy conditions for the existence and uniqueness of a solution of Eq.~(\ref{Ayd01.1u}) (see, for example,   \cite{GS_68-e}).

Following \cite{GS_68-e} we use a notation $\int  $ instead of  $\int _{{\rm {\mathbb R}}\left(\gamma \right)} $.

\begin{dfn}{\rm \cite{D_02-e,11_KchIto-e}}\label{Ayddf1}
Let a random function $u(t;x;\omega )$  and the solution of Eq.~(\ref{Ayd01.1u}) both
be defined on the same probability space. The
function $u(t; x ;\omega)$ is called a stochastic first integral for GSDE (\ref{Ayd01.1u}) if a condition
\begin{equation}\label{xGrindEQ__1_}
u\bigl(t; x (t;  x_{0});\omega\bigr)=u\bigl(0; x_{0}\bigr)
\end{equation}
is fulfilled with probability one for each solution
 $ x (t; x_{0};\omega)$  of   Eq.~(\ref{Ayd01.1u}).
\end{dfn}

Further, we shall not write a parameter $\omega$ as it is accepted.

Now we represent  the generalized It\^o-Wentzell formula which we will use for obtaining the equation for the SFI and  other results.
\begin{theorem} {\rm(the generalized It\^o-Wentzell formula)} {\rm \cite{13_KchArxDir-e,12_KchDubPrep-e}}\label{t1}
Let the random process $ x (t)\in  {\mathbb R}^{n}$ be subjected to Eq.~(\ref{Ayd01.1u}), and coefficients of this equation satisfy the conditions~${\bf \mathcal{L}_{1}}$:
\begin{description}
  \item[$ \mathcal{L}_{1}.1.$] $\displaystyle\int_{0}^{T}|a_{i} (t)| dt<\infty;$ \ \
 $\displaystyle\int _{0}^{T}|b_{i k} (t)|^{2} dt<\infty$;
  \item[$ \mathcal{L}_{1}.2.$] $\displaystyle\int_{0}^{T}dt\int   |g_{i} (t;\gamma)|^{s} \Pi (d\gamma)<\infty, \ \ s=1,2.$
\end{description}
 If $F(t; x )$, $(t; x )\in[0,T]\times\mathbb{R}^{n} $ is a scalar function which has the next stochastic differential
\begin{equation} \label{GrindEQ__2_5_1_}
d_{t} F(t; x )=Q(t; x )dt+D_{k} (t; x )dw_{k} (t)+\int G(t; x ;\gamma ) \nu (dt;d\gamma )
\end{equation}
and coefficients of Eq.~(\ref{GrindEQ__2_5_1_}) satisfy conditions ${\bf \mathcal{L}_{2}}$:
\begin{description}
  \item[$ \mathcal{L}_{2} .1.$] $Q(t; x )$, $D_{k}(t; x )$, $G(t; x ;\gamma)$ are random functions  measurable with respect to  $\mathcal{F}_{t} $, which are accorded with the processes $ w_{k}(t)$, $k=\overline{1,m}$, and  $\nu(t;\mathcal{A})$ from Eq.~(\ref{Ayd01.1u}) for each set  $\mathcal{A}\in \mathfrak{B}$ from the fixed  $\sigma-$algebra  {\rm(\cite{GS_68-e})};
  \item[$ \mathcal{L}_{2}.2.$] $
 Q (t; x )\in \mathcal{C}_{t,x}^{1,2}, \ \ D_{k}(t; x )\in \mathcal{C}_{t,x}^{1,2}, \ \ \ G(t; x ;\gamma)\in \mathcal{C}_{t,x,\gamma}^{1,2,1}
$.
\end{description}
then the next stochastic differential there exists:
\begin{equation} \label{GrindEQ__2_5_2_}
\begin{array}{c}
 \displaystyle d_{t} F(t; x (t))=Q(t; x (t))dt+D_{k} (t; x (t))dw_{k} +\\
+\displaystyle \Bigl[a_{i} (t)\frac{\partial F(t; x )}{\partial x_{i} } +\frac{1}{2} b_{ik} (t)b_{j,k} (t)\frac{\partial ^{{\kern 1pt} 2} F(t; x )}{\partial x_{i} \partial x_{j} }\Bigr. + \\
+\Bigl.\displaystyle
b_{ik} (t)\frac{\partial D_{k} (t; x )}{\partial x_{i} }\Bigr]\bigr|_{ x = x (t)} dt+
b_{ik} (t)\displaystyle\frac{\partial F(t; x )}{\partial x_{i} }\bigr|_{ x = x (t)} dw_{k} + \\
+\displaystyle\int \Bigl[(F(t; x (t)+g(t;\gamma ))-F(t; x (t))\Bigr]\nu (dt;d\gamma ) +\\
+ \displaystyle\int G\bigl(t; x (t)+g(t;\gamma );\gamma\bigr)\nu (dt;d\gamma ). \end{array}
\end{equation}
\end{theorem}

Let us remark that the generalized It\^o-Wentzell formula will be applied with extra bounding  for coefficients of Eq.~(\ref{Ayd01.1u}) \cite{11_KchOboz-e}:
\begin{equation}
 a_{i}(t)=a_{i} (t; x )\in \mathcal{C}_{t,x}^{1,2}, \ \ b_{k}(t)=b_{k}(t; x )\in \mathcal{C}_{t,x}^{1,2}, \ \ \ g_{i}(t;\gamma)=g_{i}(t; x ;\gamma)\in \mathcal{C}_{t,x,\gamma}^{1,2,1}.
\end{equation}

For simplicity, we shall use a notation $\displaystyle \frac{\partial f(t;x(t))}{\partial x_{j} }$ instead of $\displaystyle \frac{\partial f(t;x )}{\partial x_{j} }\Bigl|_{x=x(t)}$.

Let us construct an equation for the SFI of the GSDE, taking into account  that the stochastic differential for the SFI there exists and has a form:
\begin{equation} \label{xGrindEQ__2_}
d_{t} u(t;x)=Q(t;x)dt+D_{k} (t;x)dw_{k} (t)+\int G(t;x;\gamma ) \nu (dt;d\gamma )
\end{equation}

Since the definition \ref{Ayddf1} and representation (\ref{xGrindEQ__2_}) take place, we apply the generalized It\^o-Wentzell formula (\ref{GrindEQ__2_5_2_}):
\begin{equation} \label{xGrindEQ__3_}
\begin{array}{c}
  \displaystyle d_{t}u(t;x(t))=Q(t;x(t))dt+D_{k} (t;x(t))dw_{k} +b_{ik} (t)\frac{\partial }{\partial x_{i} } u(t;x(t))dw_{k} + \\
  \displaystyle  +\Bigl[a_{i} (t)\frac{\partial }{\partial x_{i} } u(t;x(t))+\frac{1}{2} b_{ik} (t)b_{j,k} (t)\frac{\partial ^{{\kern 1pt} 2} u(t;x(t))}{\partial x_{i} \partial x_{j} } + \\
  \displaystyle  +b_{ik} (t)\frac{\partial }{\partial x_{i} } D_{k} (t;x(t))\Bigr]dt+\int G(t;x(t)+g(t;\gamma );\gamma )\nu (dt;d\gamma )+\\
   \displaystyle  +\int \bigl[(u(t;x(t)+g(t;\gamma ))-u(t;x(t))\bigr]\nu (dt;d\gamma )=0.
\end{array}
\end{equation}

The  appearance of equation for the SFI depends   if the function $g(t;\gamma )$ depends or does not depend on $x$. Let us consider both situations.

 Let  function $g(t;\gamma )$ be not depend on $x$, and coefficients of Eq.~(\ref{xGrindEQ__2_}) are:
\begin{subequations}\label{xGrindEQ__4_}
\begin{align}
 \begin{array}{c}  \displaystyle{Q(t;x)=-\bigl[a_{i} (t)\frac{\partial }{\partial x_{i} } u(t;x)+\frac{1}{2} b_{ik} (t)b_{j,k} (t)\frac{\partial ^{ 2} u(t;x)}{\partial x_{i} \partial x_{j} } +} \\  \displaystyle{-b_{ik} (t)\frac{\partial }{\partial x_{i} } (b_{j,k} (t)\frac{\partial }{\partial x_{j} } u(t;x))\bigr]},
 \end{array}
 \label{xGrindEQ__4_a}\\
 \displaystyle D_{k} (t;x) = -b_{ik} (t)\frac{\partial }{\partial x_{i} } u(t;x),  \ \ \ \ \ \ \ \ \ \ \ \ \ \ \ \ \label{xGrindEQ__4_b} \\
%\label{a45}
   \displaystyle G(t;x;\gamma )=u(t;x-g(t;\gamma ))-u(x;t).\ \ \ \ \ \ \ \ \ \ \ \label{xGrindEQ__4_c}
\end{align}
\end{subequations}

\begin{theorem} In order for the random function $u(t,x)\in \mathcal{C}_{t,x}^{1,2}$ having the stochastic differential, appearing as Eq.~(\ref{xGrindEQ__2_}) to be the SFI for the system (\ref{Ayd01.1u}),  	it is sufficient that this function would be the solution of Eq.~(\ref{Ayd01.1u}) with coefficients  (\ref{xGrindEQ__4_}).
\end{theorem}

\begin{proofs} {The sufficiency} of it connects with the test of truth for Eq.~(\ref{xGrindEQ__1_}). It is fulfilled if the coefficients of Eq.~(\ref{xGrindEQ__2_}) are determined by (\ref{xGrindEQ__4_}). This fact is tested by substitution (\ref{xGrindEQ__3_}) into (\ref{xGrindEQ__4_}).~$\lozenge$
\end{proofs}

Assume that $g(t;\gamma )=g(t;x;\gamma )$. Then:
\begin{equation} \label{xGrindEQ__5_}
dx(t)=a(t)dt+b_{k} (t)dw_{k} (t)+\int g(t;x(t);\gamma ) {\rm \; }\nu (dt;d\gamma ),
\end{equation}
where $a(t)=(a_{i}(t))$, $b_{k}(t)=(b_{ik}(t))$, $g(\cdot)=(g_{i}(\cdot))$, $i=\overline{{1,n}}$,  $k=\overline{{1,m}}$.

Let us introduce new functions and make a few remarks about their properties.

If we have $y=x+g(t;x;\gamma )$ in Eq.~(\ref{xGrindEQ__3_}), then denote its solution with respect to $x$ as $x^{-1} (t;y;\gamma )$.

If we consider the domain of the uniqueness, then we can choose \begin{equation}\label{6a}
 y=z+g(t;z;\gamma ),                                                        \end{equation}
from this domain, and we get the next equality:
\begin{equation}\label{6b}
x^{-1} (t;z+g(t;z;\gamma );\gamma )=z.                                                   \end{equation}

Let   $u(t;x)$ be a random function, and its stochastic differential defined by Eq.~(\ref{xGrindEQ__2_}), and coefficients of Eq.~(\ref{xGrindEQ__2_}) are defined by  (\ref{xGrindEQ__4_a}), (\ref{xGrindEQ__4_b})  and
\begin{equation}\label{xGrindEQ__7_}
 G(t;x;\gamma )=u(t;x-g(t;x^{-1} (t;x;\gamma );\gamma ))-u(x;t).
\end{equation}

\begin{theorem} The random function $u(t;x)\in \mathcal{C}_{t,x}^{1,2}$, having the stochastic differential, appeared as Eq.~(\ref{xGrindEQ__2_}) with coefficients $D_{k} (t;x)$, $Q(t;x)$ and $G(t;x;\gamma )$ of Eq.~(\ref{xGrindEQ__2_}) defined by  (\ref{xGrindEQ__4_a}), (\ref{xGrindEQ__4_b})  and (\ref{xGrindEQ__7_}) respectively, is the  the SFI for the system (\ref{xGrindEQ__5_}). Moreover, these conditions are necessary and sufficient, if the given bounding for coefficients is set.
\end{theorem}

\begin{proofs}  Sufficiency is based on the fact, that a multiplier for the Poisson measure in Eq.~(\ref{xGrindEQ__3_}) must equal to zero. Taking into account  Eq.~(\ref{xGrindEQ__7_}) and properties   (\ref{6a}) and (\ref{6b}), we have the next equality:
\begin{equation}
\begin{array}{c}
   \displaystyle\int G(t;x(t)+g(t;x(t);\gamma ))\nu (dt;d\gamma )+ \\
  + \displaystyle \int \bigl[(u(t;x(t)+g(t;x(t);\gamma ))-u(t;x(t))\bigr]\nu (dt;d\gamma ) = \\
  = \displaystyle\int \bigl[(u(t;x(t)+g(t;x(t);\gamma )-g(t;x^{-1} (t;x+g(t;x(t);\gamma );\gamma )- \\
  -u(x(t)+g(t;x(t);\gamma );t)\bigr]\nu (dt;d\gamma )+\\
  + \displaystyle\int \bigl[(u(t;x(t)+g(t;x(t);\gamma ))-u(t;x(t))\bigr]\nu (dt;d\gamma ) \equiv 0.
\end{array}
\end{equation}

Sufficiency is proved.

Necessity follows from  conditions for  existence and uniqueness for solutions of Eq.~(\ref{xGrindEQ__5_}). Then we have  the uniqueness of the appearance for the stochastic differential $d_{t}F(t;x(t))$ for the function $F(t;x)$ which has the stochastic differential  (\ref{xGrindEQ__1_}).

Proof is complete. $\lozenge$
\end{proofs}

\section{Local stochastic density, kernel function for integral invariants for the It\^o equations and their equations}

Let $x(t)$ be a random process defined on ${\mathbb R}^{n} $ and it is a solution of the next equation
\begin{equation} \label{Ayd01.1}
\begin{array}{c}
\displaystyle{dx_{i} (t)=a_{i} (t;x(t))
dt+b_{ik} (t;x(t))
dw_{k} (t)+\int g_{i} (t;x(t);\gamma ) \nu (dt;d\gamma ),} \\
\displaystyle {x(t)=x(t;x_{0})\bigr|_{t=0} =x_{0}, \ \ \ {\mbox {\rm for every}} \ \ x_{0}\in {\rm {\mathbb R}}^{n} }, \ \ i =\overline{{1,n}}, \ \ k=\overline{{1,m}},
\end{array}
\end{equation}
and its coefficients satisfy more rigid conditions, then we used before  \cite[с.~278--290, 298--302]{GS_68}:
\begin{equation}\label{exiseq}
 a_{i}(t; x )\in \mathcal{C}_{t,x}^{1,1}, \ \ b_{ik}(t; x )\in \mathcal{C}_{t,x}^{1,2};
\end{equation}
and
\begin{equation} \label{pu}
\begin{array}{c}
\displaystyle{\int _{0}^{T}dt \int  |\nabla ^{\alpha } g(t;x;\gamma )| ^{\beta } \Pi (d\gamma )<\infty ,\ \ \ \alpha = 1,2 , \ \ \beta =\overline{1,4}, }
 \end{array}
\end{equation}
where $\nabla ^{k} $ denotes all possible varieties of combinations of $k-$th partial derivatives with respect to   $x$, taking into account a continuity of there derivatives.

Conditions (\ref{exiseq}) and (\ref{pu}) are sufficient for  existence and  uniqueness of Eq.~(\ref{Ayd01.1}) \cite{GS_68-e}.

Let us introduce a few definitions.

\begin{dfn}{\rm \cite{D-dis_79-e}}
Let $S(t)=S(t;v)$, $v\subset\hat{\mathcal{F}}\subset\mathcal{F}$ be measurable mapping from $\mathbb{R}^{n}$ into $\mathbb{R}^{n}$. The set
$
\Delta (t) = S(t;v)\cdot \Delta(0)
$
is   the dynamical invariant of the domain   $\Delta(0)$ for the process $x(t)\in \mathbb{R}^{n}$, if the following holds:
$$
\mathbf{P}\Bigl\{ x(t)\in \Delta(t)\bigl| x=x_{0}\Bigr\}=1, \ \ \ {\mbox {\rm for every}} \ \ t>0, \ \ \ {\mbox {\rm foe every}} \ \ x\in\Delta(0).
$$
\end{dfn}

Let $\rho (t;x;\omega )$ be a random function, which is ${\rm {\mathcal F}}_{t}-$measurable. Further, we shall not write a parameter $\omega$ as we said  before.

\begin{dfn}
A random function $\rho (t;x)$  is a stochastic density of the dynamic invariant for the equation connected with Eq.~(\ref{Ayd01.1}), if a random function $u(t;J;x)=J\rho (t;x)$  is the SFI,  and it satisfies the next equality
\begin{equation}\label{GrindEQ_D_2_}
J(t;x_{0})\rho (t;x(t;x_{0}))=\rho (0;x_{0}),\ \ \ {\mbox {\rm for every}}\ \ x_{0}\in \Gamma\subset R^{n},
\end{equation}
where function $J (t )=J (t,x_{0} )$ is a solution of equation
\begin{equation} \label{GrindEQ_D_3_}
\begin{array}{c}
\displaystyle{dJ (t )=J (t ) \left\{{\rm K}(t)dt+\frac{\partial b_{ik} (t;x(t) )}{\partial x_{i}} dw_{k} (t ) \right.+}\\
 +\displaystyle\left.\int \Bigl(\det \bigl[A(\delta _{i,j} +\frac{\partial g_{i} (t;x (t);\gamma )}{\partial x_{j} } )\bigr]-1\Bigr)\nu (dt,d\gamma ) \right\}, \\
\displaystyle{J (t )=J (t,x_{0} )\bigl|_{t=0} =1,{\rm \; \; dim}A(\cdot )=n\times n;}
\end{array}
\end{equation}
$\delta _{i,j} $ is  Kronecker symbol;
$$
\begin{array}{c}
K(t)=\displaystyle\Bigl[\frac{\partial a_{i}(t;x(t) )}{\partial x_{i} } + \\
\displaystyle  +\frac{1}{2} \left(\frac{\partial b_{ik} (t;x(t) )}{\partial x_{i} } \cdot \frac{\partial b_{j,k} (t;x (t) )}{\partial x_{j} } -\frac{\partial b_{ik} (t;x(t) )}{\partial x_{j} } \cdot \frac{\partial b_{ik} (t;x(t) )}{\partial x_{i} } \right) .
\end{array}
$$
and $x(t)$ is a solution of Eq.~(\ref{Ayd01.1}).
\end{dfn}

The equation (\ref{GrindEQ_D_3_}) is the one for a Jacobian of transformation  from $ x_{0}$ to $ x(t;x_{0})$. The solution of Eq.~(\ref{GrindEQ_D_3_}) can be represented in the form:
$$
\begin{array}{c}
J(t)=\displaystyle\exp \Bigl\{\int_{0}^{t}\Bigl[K(\tau ) -\frac{1}{2} \left(\frac{\partial b_{ik} (\tau ;x(\tau) )}{\partial x_{i} } \right)^{2} \Bigr]d\tau +\int_{0}^{t}\frac{\partial b_{ik}  (\tau ;x (\tau) )}{\partial x_{i} } dw_{k} (\tau ) +    \\
   +\displaystyle \int _{0}^{t}\int
   \ln \Bigl|\det \bigl[A(\delta _{i,j} +\frac{\partial g_{i} (\tau ;x(\tau)
   ;\gamma )}{\partial x_{j} } )\bigr]\Bigr|\nu (d\tau ,d\gamma ) .
\end{array}
$$

Note that $J (t )>0$ for every $t\ge 0$.

Now we define an appearance of a partial differential equation for the function  $\rho (t;x)$.
Following  Eq.~(\ref{GrindEQ_D_2_}), we get
\begin{equation} \label{GrindEQ_D_4_}
d_{t}J(t;x_{0})\rho (t;x(t;x_{0}))=0, \ \ \ {\mbox {\rm for every}} \ \ x_{0}\in \Gamma \subset  {\mathbb R} ^{n}
\end{equation}

Assume, the function $\rho (t;x)$ is a solution of this equation:
\begin{equation} \label{GrindEQ_D_5_}
\begin{array}{c}
\displaystyle{d_{t} \rho (t;x)=Q(t;x)dt+D_{k} (t;x)dw_{k} (t)+\int _{}G(t;x;\gamma ) \nu (dt;d\gamma ),} \\
\displaystyle{\rho (t;x)\left|{}_{t=0} =\right. \rho (x)\in C_{0}^{2} .}
\end{array}
\end{equation}
Here the function $\rho ( x)$ is an initial condition. The coefficients of this equation satisfy conditions
\begin{equation}\label{exiseqQQ}
 Q(t; x )\in \mathcal{C}_{t,x}^{1,1}, \ \ D_{k}(t; x )\in \mathcal{C}_{t,x}^{1,2};
\end{equation}
and
\begin{equation} \label{puGG}
\begin{array}{c}
\displaystyle{\int _{0}^{T}dt \int|\nabla ^{\alpha } G(t;x;\gamma )| ^{\beta } \Pi (d\gamma )<\infty ,\ \ \ \alpha = 1,2 , \ \ \beta =\overline{1,4}}.
 \end{array}
\end{equation}

Now we show that the function $\rho (t;x)$ would be a solution for the next equation which we obtained before in \cite{12_KchDubPrep-e} by using another approach:
\begin{equation} \label{GrindEQ_D_6_}
\begin{array}{c}
\displaystyle{d _{t} \rho (t;x)=-\left[\frac{\partial \rho (t;x)a_{i} (t;x)}{\partial x_{i} } -\frac{\partial ^{2} \rho (t;x)b_{ik} (t;x)b_{j,k} (t;x)}{\partial x_{i} \partial x_{j} } \right]dt-} \\
\displaystyle{-\frac{\partial \rho (t;x)b_{ik} (t;x)}{\partial x_{i} } dw_{k} (t)+} \\
\displaystyle{+\int \Bigl[\rho \bigl(t;x-g(t;x^{-1} (t;x;\gamma );\gamma )\bigr)\bar{D}\bigl(x^{-1} (t;x;\gamma )\bigr)-\rho (t;x)\Bigr] \nu (dt;d\gamma ),} \\
\displaystyle{\rho (t;x)\bigl|_{t=0} =  \rho (x)\in C_{0}^{2} ; \ \ \ \rho (x)\ge 0, }
\end{array}
\end{equation}
where $x^{-1} (t;x;\gamma )$ is defined as a solution of equation
\begin{equation}\label{GrindEQ_D_7_}
y+g(t;y;\gamma )=x
\end{equation}
with respect to  $y$ for every domain of uniqueness for this solution;  and $\bar{D}\left(x^{-1} (t;x;\gamma )\right)$ denotes a  Jacobian of transformation  which corresponds to  this substitution. Further, we will denote matrices' determinants as letters with bars (i.e., $\bar{A}$, $\bar{D}$).

Since the stochastic differential --  i. e. the GSDE  -- can be represented as a sum of the two paths, where the first addend is defined by the Wiener process, and the second one is defined by the Poisson jumps, hence, we rewrite  equality (\ref{GrindEQ_D_4_}) as two equalities together:
\begin{equation} \label{GrindEQ_D_8_}
\bigl[d_{t}J(t;x_{0})\rho \bigl(t;x(t;x_{0})\bigr)\bigr]_{1} =0,\ \ \ {\mbox {\rm for every}} \ \  x_{0}\in \Gamma,
\end{equation}
\begin{equation} \label{GrindEQ_D_9_}
\bigl[d_{t}J(t;x_{0})\rho \bigl(t;x(t;x_{0})\bigr)\bigr]_{2} =0, \ \ \ {\mbox {\rm for every}} \ \  x_{0}\in \Gamma ,
\end{equation}
where (\ref{GrindEQ_D_8_}) corresponds to It\^o's stochastic differential, and (\ref{GrindEQ_D_9_})  belongs to the Poisson part.

\begin{theorem} Assume that the coefficients of Eq.~(\ref{Ayd01.1}) and Eq.~(\ref{GrindEQ_D_5_}) satisfy terms (\ref{exiseq}), (\ref{pu}) and  (\ref{exiseqQQ}), (\ref{puGG}),  respectively. Then the coefficients $Q(t;x)$, $D_{k} (t;x)$, $G(t;x;\gamma )$ of Eq.~(\ref{GrindEQ_D_5_}) providing   condition (\ref{GrindEQ_D_4_}) are uniquely determined by the next equalities
\begin{description}
  \item[$I.1.$] $-Q(t;x)=\displaystyle \frac{\partial \rho (t;x)a_{i} (t;x)}{\partial x_{i} } -\frac{\partial ^{2} \rho (t;x)b_{ik} (t;x)b_{j,k} (t;x)}{\partial x_{i} \partial x_{j} };$
  \item[$I.2.$] $-D_{k} (t;x)=\displaystyle \frac{\partial \rho (t;x)b_{ik} (t;x)}{\partial x_{i} };$
  \item[$I.3.$] $G(t;x;\gamma )= \displaystyle \rho \left(t;x-g(t;x^{-1} (t;x;\gamma );\gamma )\right)\bar{D}\left(x^{-1} (t;x;\gamma )\right)-\rho (t;x).$
\end{description}
Moreover, if we have  initial conditions for Eq.~ (\ref{GrindEQ_D_6_}), then the function   $\rho (t;x)$ is the unique  solution of Eq.~(\ref{GrindEQ_D_5_}) with coefficients defined by {\rm I.1, I.2.}

\end{theorem}

\begin{proofs} The bounding of the theorem for the coefficients of Eq.~(\ref{GrindEQ_D_5_}) provides an application of the Generalized It\^o -- Wentzell formula (\ref{GrindEQ__2_5_2_}). Therefore, we get:
\begin{equation} \label{durad}
\begin{array}{c}
  \displaystyle d_{t}\rho (t;x(t))=Q(t;x(t))dt+D_{k} (t;x(t))dw_{k} +b_{ik} (t;x(t))\frac{\partial }{\partial x_{i} } \rho (t;x(t))dw_{k} + \\
  \displaystyle {+\Bigl[a_{i} (t;x(t))\frac{\partial }{\partial x_{i} } \rho (t;x(t))+\frac{1}{2} b_{ik} (t;x(t))b_{j,k} (t;x(t))\frac{\partial ^{ 2} \rho (t;x(t))}{\partial x_{i} \partial x_{j} } +} \\
 \displaystyle {+b_{ik} (t;x(t))\frac{\partial }{\partial x_{i} } D_{k} (t;x(t))\Bigr]dt+\int G\bigl(t;x(t)+g(t;x(t);\gamma )\bigr)\nu (dt;d\gamma )+}  \\
  +\displaystyle \int \Bigl[\bigl(\rho (t;x(t)+g(t;x(t);\gamma)\bigr)-\rho (t;x(t))\Bigr]\nu (dt;d\gamma).
\end{array}
\end{equation}
Let us consider the components of It\^o's differential for the above expression:
\begin{equation} \label{GrindEQ_D_10_}
\begin{array}{c}
 \bigl[d_{t}\rho (t;x(t))\bigr]_{1} =Q(t;x(t))dt+D_{k} (t;x(t))dw_{k}(t) +\\
  \displaystyle {+\bigl[a_{i} (t;x(t))\frac{\partial }{\partial x_{i} } \rho (t;x(t))+\frac{1}{2} b_{ik} (t;x(t))b_{j,k} (t;x(t))\frac{\partial ^{{\kern 1pt} 2} \rho (t;x(t))}{\partial x_{i} \partial x_{j} } +} \\
  \displaystyle {+b_{ik} (t;x(t))\frac{\partial }{\partial x_{i} } D_{k} (t;x(t))\bigr]dt+b_{ik} (t;x(t))\frac{\partial }{\partial x_{i} } \rho (t;x(t))dw_{k}(t) .}
\end{array}
\end{equation}

An appearance of  $\bigl[d_{t}\rho (t;x(t))\bigr]_{2} $ we will research later.

By using (\ref{GrindEQ_D_3_}), (\ref{GrindEQ_D_10_}) and It\^o's formula  \cite{GS_68-e}, we obtain:
$$
\begin{array}{c}
  \bigl[d_{t}J(t)\rho (t;x(t))\bigr]_{1} =\rho (t;x(t))\bigl[dJ(t)\bigr]_{1} +J(t)\bigl[d\rho (t;x(t))\bigr]_{1} + \\
  \displaystyle   +J(t)b_{ik} (t;x(t))\frac{\partial b_{j,k} (t;x(t))}{\partial x_{j} } \frac{\partial \rho (t;x(t))}{\partial x_{i} } dt= \\
  =J(t)  \displaystyle \Bigl[b_{ik} (t;x(t))\frac{\partial b_{j,k}  (t;x(t))}{\partial x_{j} } \frac{\partial \rho (t;x(t))}{\partial x_{i} } +a_{i} (t;x(t))\frac{\partial }{\partial x_{i} } \rho (t;x(t))+ \\
    \displaystyle {+\frac{1}{2} b_{ik} (t;x(t))b_{j,k} (t;x(t))\frac{\partial ^{2} \rho (t;x(t))}{\partial x_{i} \partial x_{j} } +b_{ik} (t;x(t))\frac{\partial }{\partial x_{i} } D_{k} (t;x(t))+} \\
    \displaystyle   {+b_{ik} (t;x(t))\frac{\partial }{\partial x_{i} } D_{k} (t;x(t))+Q(t;x(t))+}\\
  +  \displaystyle \rho (t;x(t))\frac{\partial a_{i} (t;x(t))}{\partial x_{i} } +\frac{1}{2} \bigl[\frac{\partial b_{ik}  (t;x(t))}{\partial x_{i} }  \frac{\partial b_{j,k}  (t;x(t))}{\partial x_{j} } -\frac{\partial b_{ik}  (t;x(t))}{\partial x_{j} }   \frac{\partial b_{ik}  (t;x(t))}{\partial x_{i} } \bigr]\Bigr]dt+\\
  +J(t)  \displaystyle \Bigl[b_{ik} (t;x(t))\frac{\partial \rho (t;x(t))}{\partial x_{i} } +\rho (t;x(t))\frac{\partial b_{ik} \left(t;x(t)\right)}{\partial x_{i} } +D_{k} (t;x(t))\Bigr]dw_{k} (t).
  \end{array}
$$
Since (\ref{GrindEQ_D_8_}) holds, we get the necessity of the theorem's conditions  I.1  and  I.2.

If we have the given restrictions, then  solutions of equations (\ref{GrindEQ_D_5_}), (\ref{GrindEQ_D_6_}) and (\ref{GrindEQ_D_10_})  exist and they are unique \cite{KR_82-e}.

Let us research the statement I.3 of this theorem.

Between the neighboring  Poisson jumps, the random process is subjected to Wiener's perturbations only, and an order of smoothness for function  $\rho (t;x(t))$ is conserved. At the time of the jumps the next functional is added:
\begin{equation}\label{Idt}
\begin{array}{c}
  \mathcal{I}( t)= \displaystyle\int_{0}^{t} \int G\bigl(\tau;x(\tau)+g(\tau;x(\tau);\gamma )\bigr)\nu (d\tau;d\gamma )+  \\
  +\displaystyle \int_{0}^{t}\int \Bigl[\bigl(\rho (\tau;x(\tau)+g(t;x(t);\gamma)\bigr)-\rho (\tau;x(\tau))\Bigr]\nu (d\tau;d\gamma)
 \end{array}
\end{equation}
and it leaves the smoothness property invariant.

Following theorems about properties for a partial differentiation equation \cite{KR_82-e}, we set, that the solution of this does exist on the next interval between jumps, and the order of its smoothness is conserved.

As (\ref{GrindEQ_D_9_}) holds, then we need to test that from the following  expressions
 \begin{equation}\label{dI2}
\begin{array}{c}
\displaystyle {\bigl[d_{t} \rho (t;x(t))\bigr]_{2} =\int \Bigl\{ (\rho (t;x(t)+g(t;x(t);\gamma ))-\rho (t;x(t))+ } \\
\displaystyle {+\Bigl[\rho (t;x(t)+g(t;x(t);\gamma )-g(t;x^{-1} (t;x(t)+g(t;x(t);\gamma );\gamma );\gamma ))\Bigr]\times } \\
\displaystyle {\times \bar{D}\bigl(x^{-1} (t;x(t)+g(t;x(t);\gamma ))\bigr)-\rho (t;x(t)+g(t;x(t);\gamma ))\Bigr\} \nu (dt;d\gamma )=}\\
 =-\displaystyle \int \Bigr[\rho (t;x(t))-\rho (t;x(t))\bar{D}\bigl(\frac{\partial x_{i}^{-1} (t;x+g(t;x;\gamma );\gamma )}{\partial (x_{j} +g_{j} (t;x;\gamma ))} \bigr)\Bigl] \nu (dt;d\gamma ),\\
\displaystyle  \bigl[dJ(t)\bigr]_{2} =J(t)\int \Bigl(\det \Bigl[A(\delta _{i,j} +\frac{\partial g_{i} (t;x(t);\gamma )}{\partial x_{j} } )\Bigr]-1\Bigr)\nu (dt,d\gamma )
\end{array}
\end{equation}
the next differential equals zero:
$$
\begin{array}{c}
  \bigl[d_{t} J(t)\rho (t;x(t))\bigr]_{2} = \\
  =\displaystyle  \int \Bigl\{ \bigr[\rho (t;x(t))+(\rho (t;x(t))\bar{D}\bigl(\frac{\partial x_{j}^{-1} (t;x+g(t;x;\gamma );\gamma )}{\partial (x_{i} +g_{i} (t;x;\gamma ))} \bigr)-\rho (t;x(t))\bigr] \times \\
  \times \displaystyle \bigl [J(t)+J(t)\bar{A}\bigl(\delta _{i,j} (t)+\frac{\partial g_{i} (t;x(t);\gamma )}{\partial x_{l} } \bigr)-J(t))\bigr])-J(t)\rho (t;x(t))\Bigr\} \nu (dt;d\gamma )= \\
  =J(t)\displaystyle  \int \Bigl[\rho (t;x(t)) \bar{D}\bigl(\frac{\partial x_{j}^{-1} (t;x(t)+g(t;x;\gamma );\gamma )}{\partial (x_{i} +g_{i} (t;x;\gamma ))} \bigr)\bar{A}\bigl(\frac{\partial (g_{i} (t;x(t);\gamma )+x_{i} )}{\partial x_{j} } \bigr)-\\
  - \displaystyle  \rho (t;x(t))\Bigr]\nu (dt;d\gamma ).
\end{array}
$$

In fact, by using  Eq.~(\ref{aGrindEQ__05_}), Eq.~(\ref{aGrindEQ__06_})  for the function $\rho (t;x(t))$,  we substitute the respective expression  to Eq.~(\ref{Idt}) in accordance with  (\ref{GrindEQ_D_7_}).

Hence, if the expression under integral sign equals zero, then the condition  I.3 is fulfilled. This is true, if $\bar{D}(\cdot )\bar{A}(\cdot )=1$.

Taking into account  the definition (\ref{GrindEQ_D_7_}) and representation $x^{-1} \bigl(t;x+g(t;x;\gamma )\bigr)=x$, we obtain:
$$
\begin{array}{c}
  \displaystyle \det \Bigl[D\bigl(\frac{\partial x_{i}^{-1} (t;x+g(t;x;\gamma );\gamma )}{\partial (x_{l} +g_{i} (t;x;\gamma ))} \bigr)A \bigl(\frac{\partial (g_{l} (t;x;\gamma )+x_{i} )}{\partial x_{j} } \bigr)\Bigr]= \\
  =\displaystyle  \det   S\bigl(\frac{\partial x_{i}^{-1} (t;x+g(t;x;\gamma );\gamma )}{\partial (x_{l} +g_{l} (t;x;\gamma ))} \frac{\partial (g_{l} (t;x;\gamma )+x_{l} )}{\partial x_{j} } \bigr)= \\
  =\displaystyle \det  S\bigl(\frac{\partial x_{i}^{-1} (t;x+g(t;x;\gamma );\gamma )}{\partial x_{j} } \bigr)=\det S\bigl(\delta _{i,j} \bigr)=1.
\end{array}
$$

It leads to the next statement: if we have domains, where one-by-one mapping between variables from equation $y=x+g(t;x;\gamma )$ there exists, then the solution of the equation from this theorem there exists and it is unique.

Proof is complete.
\end{proofs}

If a requirement for existence and uniqueness of the solution of Eq.~(\ref{Ayd01.1}) is fulfilled in the real space, then the next conditions are added \cite{D_02-e,D_89-e}:
\begin{equation}\label{ya}
\displaystyle \int _{{\rm {\mathbb R}}^{n} }f(x)\rho (t;x)d \Gamma (x)=\int _{{\rm {\mathbb R}}^{n} }f(x(t;y))\rho (y)d \Gamma(y), \ \ \ \int _{{\rm {\mathbb R}}^{n} }\rho (t;x)d \Gamma (x)=1
\end{equation}
for every continuous and bounded function $f(x)$, where a function  $\rho (x)= \rho (0;x)$ satisfies terms:
$$
\displaystyle \lim\limits_{|x|\to \infty} \frac{\partial^{k} \rho (t;x)}{\partial x_{i}^{k}}\Bigl|_{t=0}=0, \ \ \ \ k=0,1,2. \ \ \  i=\overline{1,n}.
$$
 The relation (\ref{ya}) is a definition of a stochastic kernel function for SII \cite{D_02-e,12_KchDubPrep-e}.

Let us note, that a  full collection of kernel functions  for SII  $\rho _{r} (t;x)$, $r=\overline{1,n+1}$ exists, similar to deterministic systems. As it is well known,  a  full collection of kernel functions of $n$-th order integral invariants (nII) is a collection of the kernels, if  another function which is a kernel function of  nII, would be represented as a function of  those kernels. Each kernel $\rho (t;x)$  is uniquely determined by an initial value $\rho (x)$ and the  full collection of kernels \cite{11_KchIto-e,D_89-e,Zubov_82-e}. For the stochastic processes a proof scheme of this result is similar to a proof scheme for deterministic systems.

It is necessary to stress that  existence of the density $\rho (t;x)$ connects with properties of a measure $\mu (t ;\Delta)$ ($\int _{{\rm {\mathbb R}}^{n} }\mu (t ;d(\Delta)) =const$) which is produced by some mapping. If a limit of $\mu (t ;\Delta)/\mu (\Delta )$ as $\Delta \downarrow 0$ exists there in some sense, then it is identified with the function $\rho (t;x)$. One must take into account that the function  $\rho (t;x)$   may not have a necessary smoothness which allows to construct Eq.~(\ref{GrindEQ_D_6_}).

\section{An application of the equations for the stochastic density function for obtaining the Kolmogorov equations}

Assume,  $x(t)$ is a solution of the next system (\ref{Ayd01.1}) with coefficients $a_{i} (t;x)$,\, $b_{ik} (t;x)$  and $g(t;x;\gamma )$ which are satisfied conditions (\ref{exiseq}) and (\ref{pu}).
 These conditions  provide the existence and uniqueness for solution of Eq.~(\ref{Ayd01.1}) under some terms for equation  which we consider below.

Let a random non-negative function  $\rho (t;x) $  and random process $x(t)$ be determined mutually with respect to  ${\mathcal F}_{t}$. The function  $\rho (t;x) $ is  the stochastic density function for the integral invariant which is connected with $x(t)$ as far as it satisfies conditions (\ref{ya}).

An equation for function  $\rho (t;x)$ we could obtain applying different approaches. Now we use Eq.~ (\ref{GrindEQ_D_6_}).

Further, we use the centered Poisson measure $\tilde{\nu }(\Delta t,\Delta \gamma )=\nu (\Delta t,\Delta \gamma )-\Delta t\Pi (\Delta \gamma )$ for Eq.~(\ref{GrindEQ_D_6_}):
\begin{equation} \label{aGrindEQ__04_}
\begin{array}{c}
 \displaystyle d_{t} \rho (t;x)=-\frac{\partial \rho (t;x)b_{i,k} (t;x)}{\partial x_{i} } dw_{k} (t)-\Bigl[\frac{\partial \left(\rho (t;x)a_{i} (t;x)\right)}{\partial x_{i} } - \\
  -\displaystyle\frac{1}{2} \frac{\partial ^{  2} \left(\rho (t;x)b_{i,k} (t;x)b_{j,k} (t;x)\right)}{\partial x_{i} \partial x_{j} } \Bigr]dt+  \\
  +\displaystyle\int \bigl[\rho \left(t;x-g(t;x^{-1} (t;x;\gamma );\gamma )\right)  {\bar D} \left(x^{-1} (t;x;\gamma )\right)-\rho (t;x)\bigr]\Pi (d\gamma )dt+\\
  +\displaystyle\int \bigl[\rho \left(t;x-g(t;x^{-1} (t;x;\gamma );\gamma )\right) {\bar D} \left(x^{-1} (t;x;\gamma )\right)-\rho (t;x)\bigr]\tilde{\nu }(dt;d\gamma ).
\end{array}
\end{equation}

Let us introduce a notation: $  \mathbf{M} [\rho (t;x)]=p(t;x)$.

After calculating  the mean function for Eq.~(\ref{aGrindEQ__04_})  we get that the function  $p(t;x)$ satisfies an equation
\begin{equation} \label{aGrindEQ__05_}
\begin{array}{c}
\displaystyle \frac{\partial p(t;x)}{\partial t} =-\frac{\partial \left(p(t;x)a_{i} (t;x)\right)}{\partial x_{i} } +\frac{1}{2} \frac{\partial ^{  2} \left(p(t;x)b_{i,k} (t;x)b_{j,k} (t;x)\right)}{\partial x_{i} \partial x_{j} } +  \\
+\displaystyle
 \int [p\left(t;x-g(t;x^{-1} (t;x;\gamma );\gamma )\right) {\bar D} \left(x^{-1} (t;x;\gamma )\right)-p(t;x)]\Pi (d\gamma ).                 \end{array}
\end{equation}

Equation (\ref{aGrindEQ__05_}) is   {\it Kolmogorov's equation for density function}.

Now we will define equations for a transition probability density for random processes which are solutions for equations (\ref{Ayd01.1}).

A density for transition probability for random process is a non-random function $p(t;x/s;y)$ which provides   equality
\begin{equation} \label{aGrindEQ__06_}
\displaystyle p(t;x)=\int _{{\mathbb R}^{n} } p(t;x/s;y)p(s;y) d\Gamma (y),\ \ \ t>s.
\end{equation}

From Eq.~(\ref{aGrindEQ__06_}) we conclude, that the function $p(t;x/s;y)$ does not depend on a distribution function of $p(s;y)$ for each  $s\ge 0$, hence, it is independent with respect to an arbitrary function $p(0;y)=\rho (y)$.

After the substitution (\ref{aGrindEQ__06_}) in Eq.~(\ref{aGrindEQ__05_}) and taking into account  an randomness of function $p(s;y)$ we get that term  (\ref{aGrindEQ__06_}) holds if the next equality is fulfilled
\begin{equation*}
\begin{array}{c}
\displaystyle \frac{\partial p(t;x/s;y)}{\partial t} = \ \ \ \\
=-\displaystyle \Bigl[\frac{\partial \left(p(t;x/s;y)a_{i} (t;x)\right)}{\partial x_{i} } -\frac{1}{2} \frac{\partial ^{  2} \left(p(t;x/s;y)b_{i,k} (t;x)b_{j,k} (t;x)\right)}{\partial x_{i} \partial x_{j} } \Bigr]+
\\
\displaystyle {+\int \Bigl[p\left(t;x-g(t;x^{-1} (t;x;\gamma );\gamma )/s;y\right)  {\bar D} \left(x^{-1} (t;x;\gamma )\right)-} \\ -p(t;x/s;y)\Bigr]\Pi (d\gamma ).
\end{array}
\end{equation*}

This is {\it the Kolmogorov forward equation} for the transition densities   connected with  Eq.~(\ref{Ayd01.1}).

Right now let us construct the backward equation for $p(t;x/s;y)$.

Since $p(t;x/s;y)$ is non-random function, we can rewrite Eq.~(\ref{aGrindEQ__06_}) as
$$
p(t;x)=\int _{{\mathbb R}^{n} } \mathbf{M}[ p(t;x/s;y)\rho (s;y)]d\Gamma(y).
$$

Taking into account this representation and the condition (\ref{ya}) for function $\rho (s;y)$ we go to the next equality
$$
\begin{array}{c}
\displaystyle{p(t;x)=\int _{{\mathbb R}^{n} } p(t;x/s;y) \mathbf{M}[\rho (s;y)]d\Gamma(y)%=} \\\displaystyle{
=\int _{{\mathbb R}^{n} }\mathbf{M}[ p(t;x/s;y(s;z))\rho (0;z)]d\Gamma(z),}
\end{array}
$$
where $y(s;z)$ is the solution if the It\^o equation (\ref{Ayd01.1}).

Since the left side of this equality does  not depend  on  $s$ we use the It\^o formula with respect to $s$ and obtain the following
$$
 0=\int _{{\mathbb R}^{n} }\mathbf{M}[ d_{s} p(t;x/s;y(s;z))\rho (0;z)]d\Gamma(z)=
$$
$$
=\int _{{\mathbb R}^{n} }\mathbf{M}\Bigl\{\Bigl[ \frac{\partial }{\partial s} p(t;x/s;y(s;z))+a_{j} (s;y(s;z))\frac{\partial }{\partial y_{j} } p(t;x/s;y(s;z))+
$$
$$
\displaystyle{+\frac{1}{2} b_{j,k} (s;y(s;z))b_{i,k} (s;y(s;z))\frac{\partial ^{2} }{\partial y_{j} \partial y_{i} } p(t;x/s;y(s;z))\Bigr]ds+}
$$
$$
+\displaystyle b_{j,k} (s;y(s;z))\frac{\partial }{\partial y_{j} } p(t;x/s;y(s;z))dw_{k} (s)+
$$
$$
+\displaystyle\int [p\bigl(t;x/s;y(s;z)+g(t;y(s;z);\gamma )\bigr)- p(t;x/s;y(s;z))]\nu (ds;d\gamma )\Bigr\} \rho (0;z)d\Gamma(z).
$$

Taking into account the property (\ref{ya}) again, we get the next expression:
\begin{equation*}
\begin{array}{c}
\displaystyle\int _{{\mathbb R}^{n} }\Bigl\{ \Bigl[ \frac{\partial }{\partial s} p(t;x/s;y)+a_{j} (s;y)\frac{\partial }{\partial y_{j} } p(t;x/s;y)+
\\
+\displaystyle\frac{1}{2} b_{j,k} (s;y)b_{ik} (s;y)\frac{\partial ^{2} }{\partial y_{j} \partial y_{i} } p(t;x/s;y)\Bigr]
+
\\
+\displaystyle\int \Bigl[p(t;x/s;y+g(s;y;\gamma ))-p(t;x/s;y)\Bigr]  \Pi (d\gamma )\Bigr\} p(s;y)d\Gamma(y)=0.
\end{array}
 \end{equation*}

The condition that $p(t;x/s;y)$ does not depend on $p(s;y)$  is fulfilled if the next result holds
$$
\begin{array}{c}
  \displaystyle \frac{\partial }{\partial s} p(t;x/s;y)+a_{j} (s;y)\frac{\partial }{\partial y_{j} } p(t;x/s;y)+ \\
  +\displaystyle \frac{1}{2} b_{j,k} (s;y)b_{i,k} (s;y)\frac{\partial ^{2} }{\partial y_{j} \partial y_{i} } p(t;x/s;y) + \\
  +\displaystyle \int \Bigl[p(t;x/s;y+g(s;y;\gamma ))-p(t;x/s;y)\Bigr]  \Pi (d\gamma )=0.
\end{array}
$$

This  is {\it the Kolmogorov backward equation} for the transition densities.

Let us consider Eq.~(\ref{Ayd01.1}). We used the It\^o equation with a random (non centered) Poisson measure
(\ref{Ayd01.1}).  If we can deal with the It\^o equation with a centered  Poisson measure, we can use a passage to it Eq.~(\ref{Ayd01.1}).

This approach implies a substitution from coefficients $a_{j} (\tau ;z)$ to
$$
\bar{a}_{j} (\tau ;z)=a_{j} (\tau ;z)-\displaystyle\int g(\tau ;z;\gamma ) \Pi (d\gamma )
$$
and  resulting equations are  the same as the equations which were obtained by I.I.Gihman and A.V.Skorohod \cite{GS_68-e}.

\section{Conclusion}

The concept of the  stochastic kernel function for the stochastic integral invariant allows us to obtain the well-known results such as the It\^o -- Wentzell formula \cite{D_89-e} and the Kolmogorov equations; it demonstrate the correctness of a theory based on this concept. Moreover, the method of the integral invariant gives the chance for developing the theory of the stochastic differential equations and its applications, for example,  first integrals and stochastic first integrals, and the program control with probability one (PCP1) \cite{D_78-e,D_98-e,12_KchDubPrep-e,09_ChUprW-e,11_KchUpr-e}.

%\renewcommand{\refname}{Библиографический список}
%\begin{spacing}{0.9}
%\bibliographystyle{plain}
%\bibliography{refs}

 \setlength\parsep{0pt} % чтоб уменьшить  вертик пробелы

%\end{spacing}

\end{document}